\documentclass{amsart}
\usepackage{amssymb}
\usepackage{epic,curves,mfpic,tikz}
\newtheorem{theo}{Theorem}[section]
\newtheorem{lemm}[theo]{Lemma}
\newtheorem{coro}[theo]{Corollary}
\newcommand{\cc}{\mathcal{C}}
\newcommand{\co}{\mathcal{O}}
\newcounter{fig}
\DeclareMathOperator{\Aut}{Aut}
\DeclareMathOperator{\Cay}{Cay}
\begin{document}

\address{2010 \emph{Mathematics Subjects Classification.} 20E06; (05C63; 20E08; 20F65)}
\title[Cutting up graphs revisited]{Cutting up graphs revisited -\\
a short proof of Stallings' structure theorem}
\author{B. Kr\"on}
\maketitle
Version 04.02.2010
\begin{abstract}
This is a new and short proof of the main theorem of classical structure tree theory. Namely, we show the existence of certain automorphism-invariant tree-decompositions of graphs based on the principle of removing finitely many edges. This was first done in ``Cutting up graphs'' \cite{Dunwoody1982}. The main ideas are based on the paper ``Vertex cuts'' \cite{Dunwoody2009} by M.J.~Dunwoody and the author.
We extend the theorem to a detailed combinatorial proof of J.R.~Stallings' theorem on the structure of finitely generated groups with more than one end.
\end{abstract}

\section{Historic background - Stallings' theorem, Wall's conjecture and structure trees}

The Seifert-van-Kampen Theorem (see \cite{VanKampen1933}) says that if a topological space can be decomposed into two open path connected spaces $C$ and $D$ then its fundamental group is a free product of the fundamental groups of the two subspaces with amalgamation over he fundamental group of their intersection $C\cap D$. Formally, if $x\in C\cap D$ then $\pi_1(C\cup D,x)=\pi_1(C,x)*_{\pi_1(C\cap D,x)}\pi_1(D,x)$.

In the 1960s and 1970s mathematicians began to consider groups as geometric objects themselves, and not as something that is defined by another geometric or topological object. Lie groups are considered as differentiable manifolds, and finitely generated groups are considered as Cayley graphs. Some geometric properties of finitely generated Cayley graphs can be regarded as properties of the group itself, because they do not depend on the choice of the finite set of generators. These are properties which are quasi-isometry invariants like the number of ends, growth, hyperbolicity, accessibility etc. The task of geometric group theory is to relate such geometric properties with algebraic properties of the group. In the 1980ies and 1990ies geometric group theory has become an own branch of mathematics.
 
The counter part of the Seifert-van-Kampen Theorem in geometric group theory is Stallings' structure theorem. A group $G$ is said to \emph{split} over a subgroup $H$ if $G$ is a non-trivial free product with amalgamation over $H$ or $G$ is an HNN-extension of a group over $H$. A finitely generated group is said to have more than one end if its finitely generated Cayley graphs have more than on end. Equivalently, if there is a finite subgraph of the Cayley graph whose complement has at least two infinite components.

\begin{theo}[Stallings' structure theorem]
A finitely generated group has more than one end if and only if it splits over some finite subgroup.
\end{theo}

Stallings has proved this theorem for the torsion free case in 1968 \cite{Stallings1968} and for the general case in 1972 \cite{Stallings1971}. For his work he was awarded the Cole prize in 1970.

It may happen that the factors of this group decomposition split again over some finite subgroup, and the factors of this second splitting may split again, and so on. A group is said to be \emph{accessible} if this process of splitting over finite subgroups stops after finitely many steps.
In 1971 C.T.C.~Wall conjectured in \cite{Wall1971} that all finitely generated groups are accessible. A first progress in solving the problem was made in 1976 by Bamford and Dunwoody \cite{Bamford1976} by finding a criterion for accessibility.

A tree which corresponds to an automorphism invariant tree-decomposition of a graph $X=(VX,EX)$ is called a \emph{structure tree}. In the present context, a \emph{structure cut} is a component $C$ in the complement of a finite set of edges such that $C$ and $VX\setminus C$ both contain a ray (one-way infinite path), and $C$ is nested with $g(C)$, for any automorphism $g$. Being nested means that $C\subset g(C)$ or $g(C)\subset C$ or $C\cap g(C)=\emptyset$ or $C\cup g(C)=VX$. Such structure cuts yield structure trees, see Sections~\ref{sec:blocks} and \ref{sec:trees}.

Structure trees were introduced in 1979 \cite{Dunwoody1979}. Three years later in the paper ``Cutting up graphs'' \cite{Dunwoody1982} the existence of structure cuts was finally proved for all graphs with more than one edge-end. These are graphs with a finite set of edges whose complement has two components each of which contain a ray.
 
When we consider the action of finitely generated groups with more than one end on structure trees of their Cayley-graphs, then Bass-Serre Theory implies Stallings' Theorem, see Section~\ref{sec:trees}.

There are also applications of structure trees in graph theory, see \cite{Hamann2009,Kroen2009,Kroen2001,Moeller1992ends1,Moeller1996,Moeller1995,Seifter2008,Thomassen1993}.
Structure cuts have been further developed by Dicks and Dunwoody in 1989 in the book \cite{Dicks1989}.
In 1985 \cite{Dunwoody1985} Wall's conjecture was proved for finitely presented groups.

In 1993 Thomassen and Woess introduced a graph theoretic notion of accessibility in \cite{Thomassen1993}. They called a graph \emph{accessible} if there is an integer $n$ such that any two ends can be separated by removing $n$ (or fewer) edges. They showed that a finitely generated group is accessible if and only if its finitely generated Cayley graphs are accessible. Dicks and Dunwoody have shown in \cite{Dicks1989} that for all $n$ there are systems of structure cuts which separate any pair of ends that can be separated by $n$ or fewer edges. Hence for Cayley graphs of accessible groups, there are structure trees which describe all possible splittings of the group.

In the same year, Wall's conjecture was finally disproved, see
~\cite{Dunwoody1993}.\smallskip

In \cite{Dunwoody2009} Dunwoody and the author have proved the existence of structure cuts which are based on the principle of removing finite sets of vertices instead of edges. These results imply a generalization of Stallings' theorem from finitely generated to arbitrarily generated groups. Namely, a group splits over a finite subgroup if and only if it has a Cayley graph with more than one vertex end. That is, a Cayley graph with two rays which can be separated from each other by removing finitely many vertices. Another application of \cite{Dunwoody2009} is the generalization of Tutte's tree decomposition of 2-connected graphs to $k$-connected graphs for any integer $k$. The arguments in \cite{Dunwoody2009} yield a proof of the classical result on the existence of structure cuts which is contained in Lemmas~\ref{lemm:finitely}, \ref{lemm:intersection}, \ref{lemm:not_nested_corner}, \ref{lemm:corners_equality} and Theorem~\ref{theo:optimally}. Together with a certain tree-construction in Section~\ref{sec:blocks} and some Bass-Serre Theory in Section~\ref{sec:trees} we obtain a complete proof of Stallings' Theorem in Section~\ref{sect:stallings}.

This way of proving Stallings' theorem is in principle not new and was also mentioned as application in \cite{Dunwoody1982}. What is new are the arguments that follow from the results in \cite{Dunwoody2009} and give a short proof of the existence of structure cuts, for instance see Lemma~\ref{lemm:not_nested_corner} and Theorem~\ref{lemm:corners_equality}. The short proof of Thomassen and Woess of Lemma~\ref{lemm:finitely} also simplified the original proof. An improvement of structure tree theory is, that the vertices of the tree are not defined as equivalence classes of cuts, but as certain inseparable blocks which are subsets of the underlying graph. This approach does not shorten the construction of structure trees significantly, but we think that it is more accessible to inexperienced readers. And last but not least, one goal of the paper is to present a complete and detailed combinatorial proof of Stallings' theorem.

\section{Minimal edge cuts}

Let $X=(VX,EX)$ be an undirected simple graph. That is, edges are two-element sets of vertices. For $C,D\subset VX$ let $\delta(C,D)$ denote the set of edges with one vertex in $C$ and one vertex in $D$. We write $C^\mathrm{c}=VX\setminus C$ and call $\delta C=\delta(C,C^\mathrm{c})$ the \emph{edge boundary} of $C$. A set of vertices $C$ is \emph{connected} if the subgraph spanned by $C$ is connected. A \emph{$k$-separator} is a $k$-element edge boundary of a set of vertices $C$, where $C$ and $C^\mathrm{c}$ are connected. For a set of edges $F$ define $X-F$ as the graph $(VX,EX\setminus F)$.

\begin{lemm}[Proposition 4.1 in \cite{Thomassen1993}]\label{lemm:finitely}
Let $e$ be an edge of a connected graph $X$ and let $k$ be an integer. There are only finitely many $k$-separators which contain $e$.
\end{lemm}

\begin{proof}
We prove the statement by induction on $k$. The case $k=1$ is obvious.

Suppose the statement holds for all connected graphs for some integer $k\ge 1$. We show the statement in $X$ for $(k+1)$-separators containing $e=\{x,y\}$. The graph $X-\{e\}$ is connected, because $k\ge 2$. Hence there is a path $\pi$ from $x$ to $y$ in $X-\{e\}$. Every $(k+1)$-separator in $X$ which contains $e$ also contains an edge $e'$ of $\pi$. By the induction hypothesis there are only finitely many $k$-separators in $X-\{e\}$ which contain $e'$. Now the statement follows, because $\pi$ is finite and different $(k+1)$-separators in $X$ which contain $e$ and $e'$ correspond to different $k$-separators in $X-\{e\}$ which contain $e'$.
\end{proof}

The \emph{boundary} $NC$ of a set of vertices $C$ is the set of vertices in $C^\mathrm{c}$ which are adjacent to some vertex in $C$. We write $\beta C$ to denote the set $NC\cup NC^\mathrm{c}=\bigcup\delta C$. 

A \emph{component} of a set of vertices $A$ is a maximal connected subset of $A$. Vertices $x,y$ are separated by $S\subset VX$ if $x,y$ lie in different components of $VX\setminus S$.

Sets of vertices $A,B$ are \emph{separated} by a set of vertices $S$  if any $x\in A$ and $y\in B$ lie in distinct components of $VX\setminus S$.
Sets of vertices $A,B$ are \emph{separated} by a set of  edges $F$ if any $x\in A$ and $y\in B$ lie in distinct components of the graph $X-F$, respectively. A vertex $x$ is said to be separated from a set of vertices $A$ (or a vertex) if $\{x\}$ is separated from $A$.

A ray is a one-way infinite path (of distinct vertices). A \emph{tail} of a ray is an infinite subpath of a ray. Two rays are said to be separated by a set (of vertices or edges) if this set separates some tails of the rays.
We call two rays \emph{edge equivalent} if they cannot be separated by a finite set of edges. The corresponding ends are called the \emph{edge ends}. In non-locally finite graphs there are different notion of ends (usually defined by separation by removing finite sets of vertices), but for locally finite graphs, all definitions coincide and they correspond to Freudenthal's end compactification of locally compact Hausdorff space \cite{Freudenthal1931,Freudenthal1942,Freudenthal1944}.

A \emph{cut} (or \emph{edge-cut}) is a set of vertices $C$ with finite edge boundary such that $C$ and $C^\mathrm{c}$ are both connected and contain a ray.
If a cut contains a ray $R$ then it contains all rays which are equivalent to $R$. Hence if $R$ lies in $C$ then we say that $C$ contains the corresponding end. 
If there is an edge cut then let $\kappa$ be the minimal cardinality of all boundaries of edge cuts. Edge cuts $C$ with $|\delta C|=\kappa$ are called \emph{minimal}. If $X$ is connected and has more than one edge end then there is a minimal edge cut.

\begin{lemm}\label{lemm:intersection}
Let $C$ and $D$ be minimal cuts. If $C\cap D$ and $C^\mathrm{c}\cap D^\mathrm{c}$ are cuts then they are minimal cuts.
\end{lemm}

In \cite[Theorem 2]{Jung1977} and \cite[Proposition 2.1]{Jung1993} Jung and Watkins prove a similar result.

\begin{proof}
According to Figure~\ref{fig:1block} we set\smallskip

\begin{tabular}{lll}
$a=|\delta(C\cap D,C^\mathrm{c}\cap D)|$,\quad& $b=|\delta(C\cap D,C\cap D^\mathrm{c})|$,& $c=|\delta(C\cap D^\mathrm{c},C^\mathrm{c}\cap D^\mathrm{c})|$,\\
$d=|\delta(C^\mathrm{c}\cap D,C^\mathrm{c}\cap D^\mathrm{c})|$,& $e=|\delta(C\cap D,C^\mathrm{c}\cap D^\mathrm{c})|$,&$f=|\delta(C\cap D^\mathrm{c},C^\mathrm{c}\cap D)|$.
\end{tabular}\smallskip

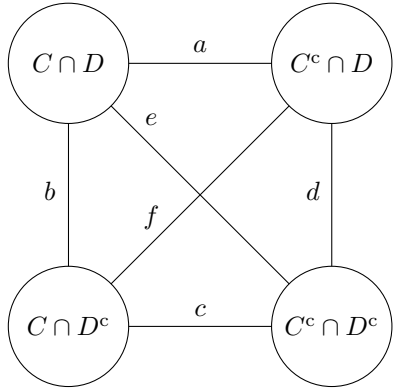
\begin{figure}[htbp]
\centering
\begin{tikzpicture}
\path (0,0) coordinate (p1);
\path (3.5,0) coordinate (p2);
\path (0,3.5) coordinate (p4);
\path (3.5,3.5) coordinate (p3);
\draw (0.8,0) -- (2.7,0);
\draw (0.8,3.5) -- (2.7,3.5);
\draw (0,0.8) -- (0,2.7);
\draw (3.5,0.8) -- (3.5,2.7);
\draw (0.565685425,0.565685425) -- (2.934314575,2.934314575);
\draw (2.934314575,0.565685425) -- (0.565685425,2.934314575);
\draw (p1) circle (0.8cm) node {$C\cap D^\mathrm{c}$};
\draw (p2) circle (0.8cm) node {$C^\mathrm{c}\cap D^\mathrm{c}$};
\draw (p3) circle (0.8cm) node {$C^\mathrm{c}\cap D$};
\draw (p4) circle (0.8cm) node {$C\cap D$};
\draw (1.75,3.72) node {$a$};
\draw (1.75,0.22) node {$c$};
\draw (-0.25,1.8) node {$b$};
\draw (3.25,1.8) node {$d$};
\draw (1.1,2.75) node {$e$};
\draw (1.1,1.45) node {$f$};

\end{tikzpicture}
\caption{One-connected graph and structure tree}\label{fig:1block}
\end{figure}
Then
\[\kappa = |\delta C|= a + e+f +c = |\delta D|= b + e + f + d\]
and hence
\begin{equation}\label{equa:4n}
2\kappa = a +b +c +d +2e+2f.
\end{equation}
The sets $C\cap D$ and $C^\mathrm{c}\cap D^\mathrm{c}$ contain an end and so
\[|\delta (C\cap D)|=a+e+b\ge\kappa\quad\mbox{and}\quad|\delta (C^\mathrm{c}\cap D^\mathrm{c})|=c+e+d\ge\kappa.\]
Hence $a +b +c +d +2e\ge 2\kappa$ and, by (\ref{equa:4n}), $a +b +c +d +2e= 2\kappa$
and $f=0$. Finally, $a+e+b=c+e+d=\kappa$ and $|\delta (C\cap D)|=|\delta (C^\mathrm{c}\cap D^\mathrm{c})|=\kappa.$
\end{proof}

\section{Main Theorem}\label{sect:maintheo}

Sets of vertices $C$ and $D$ are \emph{nested} if $C\subset D$, $C^\mathrm{c}\subset D$, $C\subset D^\mathrm{c}$ or $C^\mathrm{c}\subset D^\mathrm{c}$. Equivalently, if $C\subset D$, $D\subset C$, $C\cap D=\emptyset$ or $C\cup D=VX$. Equivalently, if one of following intersections is empty $C\cap D,C\cap D^\mathrm{c},C^\mathrm{c}\cap D,C^\mathrm{c}\cap D^\mathrm{c}$. These intersections are called \emph{corners} of $C$ and $D$. According to Figure~\ref{fig:1block}, we say that $C\cap D$ is \emph{opposite} to $C^\mathrm{c}\cap D^\mathrm{c}$, and $C^\mathrm{c}\cap D$ is \emph{opposite} to $C\cap D^\mathrm{c}$.

\begin{lemm}\label{lemm:not_nested_corner}
Let $C,D,E$ be sets of vertices and let $C$ and $D$ be not nested.

If $E$ is not nested with two opposite corners of $C$ and $D$ then $E$ is not nested with both $C$ and $D$.
If $E$ is not nested with some corner of $C$ and $D$ then $E$ is either not nested with $C$ or not nested with $D$.
\end{lemm}

\begin{proof}
Suppose $E$ is not nested with two opposite corners.
By relabeling $C^\mathrm{c}$ as  $C$, if necessary, we  can assume that
$E$ is not nested with $C\cap D$ and $C^\mathrm{c}\cap D^\mathrm{c}$.

Suppose $E$ and $C$ are nested. If $E\subset C$ or $E^\mathrm{c}\subset C$ then this would contradict the assumption that $E$ is not nested with $C^\mathrm{c}\cap D^\mathrm{c}$. If $E\subset C^\mathrm{c}$ or $E^\mathrm{c}\subset C^\mathrm{c}$ then this would contradict the assumption that $E$ is not nested with $C\cap D$.
The assumption that $E$ and $D$ are nested leads to a contradiction in the same way. Hence the first claim is established.

Now suppose that $E$ is not nested with some corner, say $C\cap D$, and suppose $E$ is nested with both $C$ and $D$. There are four possible inclusions for $E$ and $C$ being nested, and four for $E$ and $D$ being nested. We show that the corresponding 16 possibilities all lead to a contradiction. If $C\subset E$ or $D\subset E$ then  $C\cap D\subset E$, contradicting $E$ not being nested with $C\cap D$. Also if $C\subset E^\mathrm{c}$ or $D\subset E^\mathrm{c}$ then  $C\cap D\subset E^\mathrm{c}$. Now 4 cases remain. If one of the sets $C^\mathrm{c}$ and $D^\mathrm{c}$ were in $E$ and the other in $E^\mathrm{c}$, then $C^\mathrm{c}\cap D^\mathrm{c}=\emptyset$, and $C$ and $D$ would be nested. 
If $C^\mathrm{c}$ and $D^\mathrm{c}$ are both in $E$ or both in $E^\mathrm{c}$ then $C^\mathrm{c}\cup D^\mathrm{c}\subset E$ or $C^\mathrm{c}\cup D^\mathrm{c}\subset E^\mathrm{c}$, implying $E^\mathrm{c}\subset C\cap D$ or $E\subset C\cap D$, respectively,
and $C\cap D$ would be nested with $E$.
\end{proof}

Let $C$ be a cut and let $M(C)$ be the set of minimal cuts which are not  nested with $C$. Set $m(C) = |M(C)|$. It follows from Lemma~\ref{lemm:finitely}  that $m(C)$ is finite.

\begin{lemm}\label{lemm:corners_equality}
Let $C$ and $D$ be edge-cuts which are not  nested and suppose $C\cap D$ and $C^\mathrm{c}\cap D^\mathrm{c}$ are cuts, then\[m(C\cap D) + m(C^\mathrm{c}\cap D^\mathrm{c}) < m(C) + m(D).\]
\end{lemm}

\begin{proof}
It follows from Lemma~\ref{lemm:intersection} that $C\cap D$ and $C^\mathrm{c}\cap D^\mathrm{c}$ are minimal cuts.

Let $E$ be a minimal cut. If $E$ is in $M(C^\mathrm{c}\cap D^\mathrm{c})\cap M(C\cap D)$ then, by Lemma~\ref{lemm:not_nested_corner}, $E$ is in $M(C)$ and in $M(D)$. Hence if $E$ is counted twice on the left of the above inequality then it is also counted twice on the right.

If $E$ is in $M(C\cap D)\setminus M(C^\mathrm{c}\cap D^\mathrm{c})$ or in $M(C^\mathrm{c}\cap D^\mathrm{c})\setminus M(C\cap D)$, that is $E$ is counted once on the left, then, again by Lemma~\ref{lemm:not_nested_corner}, $E$ is in $M(C)$ or in $M(D)$. Hence $E$ is counted at least once on the right side of the inequality. We have now proved that $m(C\cap D) + m(C^\mathrm{c}\cap D^\mathrm{c}) \le m(C) + m(D)$. Since $C\in M(D)$ and $D\in M(C)$, the cuts $C$ and $D$ are counted on the right side, but not on the left side, so this inequality is a strict inequality.
\end{proof}

Let $\cc$ be the set of all minimal cuts. Set $m=\min\{m(C)\mid C\in \cc\}$. This minimum exists as integer, because the values $m(C)$ are all finite. A minimal cut $C$ with $m(C)=m$ is called \emph{optimally nested}. The following is the main theorem in classical structure tree theory.

\begin{theo}\label{theo:optimally}
Optimally nested cuts are nested with all other optimally nested cuts.
\end{theo}

\begin{proof}
We show that optimally nested cuts are nested with all other cuts.
Suppose there are optimally nested minimal cuts $E$ and $F$ which are not  nested. Then $m\ge 1$. There cannot be two adjacent corners which do not contain an end, because $C$, $C^\mathrm{c}$, $D$ and $D^\mathrm{c}$ all contain an end. Hence there is a pair of opposite corners which contain an end. By relabeling we can assume that these corners are $C\cap D$ and $C^\mathrm{c}\cap D^\mathrm{c}$ and by 
Lemma~\ref{lemm:intersection}, each of $C\cap D$ and $C^\mathrm{c}\cap D^\mathrm{c}$ are minimal edge-cuts.
Now Lemma \ref{lemm:corners_equality} says that 

\[m(C\cap D) + m(C^\mathrm{c}\cap D^\mathrm{c})< m(C) + m(D)  = 2m. \]
Thus one of the summands on the left side is less than $m$, contradicting the minimality of $m$.
\end{proof}

\section{Blocks and trees from nested systems}\label{sec:blocks}

Let $\cc$ be a set of sets of vertices. A nonempty set of vertices $B$ is called \emph{$\cc$-inseparable} if no pair of vertices in $B$ can be separated by $\beta C$, for any $C\in \cc$. In other words, for all $C\in\cc$ either $B\subset C\cup NC$ or $B\subset C^\mathrm{c} \cup NC^\mathrm{c}$. Maximal $\cc$-inseparable sets are called the \emph{$\cc$-blocks}. Note that edges are $\cc$-inseparable and distinct blocks are not necessarily disjoint. For a block $B$, let $\cc(B)$ denote the set of all $C$ in $\cc$ which are minimal with respect to the inclusion $B\subset C\cup NC$. That is, $C$ is in $\cc(B)$ if $B\subset D\cup ND\subset C\cup NC$, for $C,D\in\cc$, implies $C=D$.

We call $\cc$ nested if any two sets in $\cc$ are nested.

\begin{lemm}\label{lemm:blocks}
Let $\cc$ be a nested set of sets of vertices and $C\in\cc$. No pair of vertices in $\beta C$ is separated by $\beta D$, for any $D\in\cc$.
Let $\cc$ be minimal.
 There is precisely one $\cc$-block $B_C$ such that $C\in\cc(B_C)$. If $D\in\cc(B_C)$ then $\beta D\varsubsetneq B_C$. Moreover,
\begin{equation}\label{equa:block}
\bigcup_{D\in\cc(B_C)}\!\beta D\ \subset\ B_C\  =\ \bigcap_{D\in\cc(B_C)}\!D\cup ND.
\end{equation}
\end{lemm}

\begin{proof}
Suppose $x,y\in\beta C$ are separated by $\beta D$. After possibly replacing $C$ with $C^\mathrm{c}$ and $D$ with $D^\mathrm{c}$, we have $C\cap D=\emptyset$ and $x\in C^\mathrm{c} \cap D$. Since $x\in\beta C$, $x$ is adjacent to some vertex in $C\cap D^\mathrm{c}$. Hence $x\in ND^\mathrm{c}$, contradicting the assumption that $\beta D$ separates $x$ from another vertex. 

Suppose there are different blocks $B,B'$ such that $C\in\cc(B)\cap\cc(B')$. Suppose there is a vertex $x\in B'\setminus B$ and $y\in B\cap B'$. They are separated by $\beta D$, for some $D\in\cc$. Any path $\pi\subset C$ from $x$ to $y$ intersects $D$ and $D^\mathrm{c}$. Hence $C\cap D\ne\emptyset$ and $C\cap D^\mathrm{c}\ne\emptyset$. So either $C^\mathrm{c}\cap D=\emptyset$ or $C^\mathrm{c}\cap D^\mathrm{c}=\emptyset$, equivalently $D\subset C$ or $D^\mathrm{c}\subset C$. One of the sets $D\cup\beta D,D^\mathrm{c}\cup\beta D$ contains $B$, the other $B'$. If $D\subset C$ and $B\subset D\cup\beta D$ then $B\subset D\cup ND\varsubsetneq C\cup NC$ in contradiction to $C\in\cc(B)$. Any of the other cases leads to a contradiction in the same way.

The intersection in (\ref{equa:block}) is maximal inseparable, it contains $NC$ and it is contained in $C\cup NC$. Hence this is the unique block $B_C$ such that $C\in\cc(B_C)$. 

If $D,E\in\cc(B_C)$ then $E^\mathrm{c}\subset D$ which implies $E^\mathrm{c}\cup\beta E\subset D\cup\beta D$ and
\[\beta E\subset \bigcap_{D\in\cc(B_C)}\!D\cup ND=B_C\]
and establishes the inclusion in (\ref{equa:block}). If $\cc(B_C)=\{C\}$ then $B_C=C\cup NC$, and hence $\beta C$ is a proper subset of $B_C$.
 If $\cc(B_C)$ contains a cut $D$, $D\ne C$, then $\beta C\cup\beta D\subset B_C$ and again $\beta D$ is a proper subset of $B_C$.
\end{proof}

Given a nested set $\cc$ of minimal cuts we define a graph $T=T(\cc)$. Let $VT$ be the set of $\cc$-blocks. Two vertices (blocks) $v,w$ of $T$ are defined to be adjacent if they intersect.

\begin{theo}
Let $\cc$ be a nested minimal system of edge cuts. Then $T(\cc)$ is a tree.
\end{theo}

\begin{proof}
Given an edge $\{v,w\}\in ET$ there is a cut $C\in\cc$, such that $\beta C\cap v\cap w\ne\emptyset$. Lemma~\ref{lemm:blocks} implies $\beta C\subset v\cap w$. The graph $T-\{v,w\}$ is disconnected and hence $T$ does not contain any circuits.

Let $B_1,B_2$ be two $\cc$-blocks. Lemma~\ref{lemm:finitely} says that any pair of vertices in these blocks can be separated in $X$ only by finitely many sets $\beta C$, for $C\in\cc$. This implies that there is a finite path from $B_1$ to $B_2$.

Let $\pi$ be a path between vertices $x,y\in VT$. Lemma~\ref{lemm:finitely} implies that there are only finitely many sets $\beta C$, $C\in\cc$, which contain some edge in $\pi$. Hence there is a finite path in $T$ connecting the blocks which contain $x$ and $y$, respectively. This means that $T$ is connected and thus $T$ is a tree.
\end{proof}

\section{Group splitting and Bass-Serre Theory}\label{sec:trees}

Let $H,J$ be groups and $A<H$, $B<J$ be isomorphic subgroups. The amalgamated product with isomorphism $ \varphi:A\to B$ is
\[H*_{A}J=\left<H,J\mid a= \varphi(a), a\in A\right>.\]
Let $T_H$ be a system of representatives of the left cosets of $A$ in $H$ and $T_J$ of left cosets of $B$ in $J$, where $A$ and $B$ are represented by the neutral element $1$. A \emph{normal form for $H*_{A}J$} is a sequence $(x_0,x_1,\ldots,x_n,a)$ such that $a\in A$ and $x_i\in T_H\setminus\{1\}\cup T_J\setminus\{1\}$, and no consecutive elements $x_i$ and $x_{i+1}$ lie in the same system of representatives.

Let $A,B$ be isomorphic subgroups of $H$. The HNN-extension with isomorphism $ \varphi:A\to B$ is
\[H*^{A}=\left<H,t\mid tat^{-1}= \varphi(a), a\in A\right>,\]
where $t$ is an additional generator, called the \emph{stable letter} which is not contained in $H$. A \emph{normal form for $H*^{A}$} is a sequence  $(x_0,t^{\varepsilon_0},x_1,t^{\varepsilon_1},\ldots ,x_n,t^{\varepsilon_n},h)$ where $h$ is an arbitrary element of $H$, $\varepsilon_i\in\{-1,1\}$, there is no consecutive subsequence $t^\varepsilon,1,t^{-\varepsilon}$ and if $\varepsilon_i=1$ then $x_i\in B$, if $\varepsilon_i=-1$ then $x_i\in A$.

Note that the notations $H*_{A}J$ and $H*^{A}$ are ambiguous, because the amalgamated product and the HNN-extension are not determined by $H,J,A,B$, they depend on the choice of $ \varphi$. The following can for instance be found as Theorems~11.3 and 14.3 in Bogopolski's book \cite{Bogopolski2008} in terms of right co-sets instead of left-cosets.

\begin{lemm}\label{lem:normalform}
For every element $g$ in a free product with amalgamation or in an HNN-extension there is a unique normal form $(a_1,a_2,\ldots,a_n)$ such that $g=a_1a_2\ldots a_n$.
\end{lemm}

\begin{proof}
Any $h\in H$ can be uniquely written as $[[x]]^H\cdot [x]^H$ where $[[x]]^H\in T_H$ and $[x]^H\in A$.
Let $W$ be the set of normal forms for $G=H*_{A}J$. 
We define an action of $H$ on $W$ on the right by $(x_0,x_1,\ldots,x_n,a)\cdot h$
\[=\begin{cases}
(x_0,x_1,\ldots,x_n,ah)&\text{if } h\in A,\\
(x_0,x_1,\ldots,x_n,[[ah]]^H,[ah]^H)&\text{if }h\notin A,\ x_n\in T_J,\\
(x_0,x_1,\ldots, x_n ah)&\text{if } h\notin A,\ x_n\in T_H,\ x_n ah\in A,\\
(x_0,x_1,\ldots, [[x_n ah]]^H,[x_n ah]^H)&\text{if } g\notin A,\ x_n\in T_H,\ x_n ah\notin A\\
\end{cases}\]
and we define
\[(a)\cdot h=\begin{cases}
(ah)&\text{if } h\in A,\\
([[ah]]^H,[ah]^H)&\text{if }h\notin A.
\end{cases}\]
We can do the same for $J$. 
The actions of $H$ and $J$ on $W$ can be extended to an action of the free product $H*J$ on $W$. In this free product, elements $a\in A$ and $\varphi(a)\in B$ are not identified, but elements of the form  $a\varphi(a)^{-1}$ are in the kernel of this action. The same holds for the normal closure $N$ of these elements. Hence we obtain a well defined action of $G=H*_{A}J=(H*J)/N$ on $W$. If an element $g\in G$ would have two different normal forms $(x_0,x_1,\ldots,x_n,a)$ and $(y_0,y_1,\ldots,y_n,a')$ then
\[(1)\cdot g=(1)\cdot x_0x_1\ldots x_na=(x_0,x_1,\ldots,x_n,a)\text{\quad  and}\]
\[(1)\cdot g=(1)\cdot y_0 y_1\ldots y_na'=(y_0,y_1,\ldots,y_n,a'),\]
which is impossible because $(1)\cdot g$ is well defined.

In the case of an HNN-extension we first define actions of $H$ and $\{t\}$ on the set of all normal forms $W$, similar to the case of amalgamated products. We obtain an action of $H*\left<t\right>$ on $W$. The action of $tat^{-1}$ and $\varphi(a)$ on $W$ coincide for $a\in A$. Hence $tat^{-1}\varphi(a)^{-1}$ is in the kernel of this action, and so this also holds for the normal closure $N$ of all such elements. Since $G=H*^A=(H*\left<t\right>)/N$, we get an action $G$ on $W$ and proceed as before.
\end{proof}

With $\Aut(X)$ we denote the automorphism group of $X$. A group $G$ is said to \emph{act} on a graph $X$ if there is a homeomorphism $\psi:G\to\Aut(X)$. We usually write $g$ instead of $\psi(g)$ if this does not cause any confusion. An action on $X$ is said to be \emph{transitive} if it is transitive on $VX$. That is for all $x,y\in VX$ there is a $g\in G$ such that $g(x)=y$. If $G$ acts on $X$ and if $\cc$ is a $G$-invariant nested set of cuts then the action of $G$ on on $X$ induces an action of $G$ on the set of blocks and hence on the tree $T(\cc)$.
If the action is transitive on $X$ then it is also transitive on $T(\cc)$.

If $G$ acts on $X$ then the vertices of the quotient graph $X/G$ are the orbits of the action of $G$ on $VX$. Two vertices $v,w$ of $X/G$ are adjacent if there are vertices $x\in v$, $y\in w$ such that $x,y$ are adjacent in $X$. We consider this quotient graph $Y$ as multigraph.
 That is, $VY$ and $EY$ are arbitrary sets and there are functions $\alpha:EY\to VY$ and $\omega:EY\to VY$ which termine origin and terminal vertex of an edge. 
A \emph{loop} is a multigraph with one vertex $x$ and one edge $e$. That is, $\alpha(e)=\omega(e)=x$. A \emph{segment} is a connected multigraph with two vertices $x,y$ and one edge $e$. That is, $\alpha(e)=x$ and $\omega (e)=y$.

An \emph{edge inversion} is an element $g\in G$ together with an edge $\{x,y\}\in EX$ such that $g(x)=y$ and $g(y)=x$.

The following theorem from Bass-Serre Theory can also be found in the books \cite{Bogopolski2008,Serre1980,Serre2003}.

\begin{theo}\label{theo:serre}
Let $G$ act without edge inversion and transitively on an infinite tree. Then $G$ splits over the stabilizer of an edge of the tree.

If $X/G$ is a segment then $G$ splits as a non-trivial free product with amalgamation over the stabilizer of an edge of $T$.
If $X/G$ is a loop then $G$ splits as HNN-extension over the stabilizer of an edge of $T$. The stable letter maps the origin vertex of that edge to the terminal vertex.
\end{theo}

Let $G_P,G_Q$ be the stabilizers of the vertices $P,Q$ and let $G_e=G_P\cap G_Q$ denote the pointwise stabilizer of the edge $e=\{P,Q\}$.

\begin{proof}
Suppose $X/G$ is a segment. Let $e=\{P,Q\}$ be an edge of $T$ and set $G'=G_P*_{G_e} G_Q$. Because $G_P\cup G_Q$ generates both $G$ and $G'$, and $G_P\cap G_Q=G_e$ in both groups $G,G'$ and because any relation in $G'$ is a relation in $G$, there is a unique homomorphism $\psi: G'\to G$ which is the identity on $G_P\cup G_Q$. The homeomorphism yields an action of $G'$ on $T$ and this action is transitive.

We have to show that $\psi$ is bijective. Surjectivity follows from the fact that $G_P\cup G_Q$ generates both groups and that any relation in $G'$ is a relation in $G$.

To see that $\psi$ is injective, choose an element $g\in G'\setminus \{1\}$. If $g\in G_P\cup G_Q$ then $g$ is not in the kernel of $\psi$ because $\psi$ is the identity on $G_P\cup G_Q$. Otherwise, if $g\in G\setminus (G_P\cup G_Q)$ then let $(x_0,x_1,\ldots,x_n,a)$ be the standard presentation of $g$ with respect to $T_1$ as a set of representatives of left cosets of $G_e$ in $G_P$ and $T_2$ as a set of representatives of left cosets of $G_e$ in $G_Q$. Suppose $x_n\in T_1\setminus\{1\}$. Then $x_n,x_{n-2},x_{n-4},\ldots$ act as rotations of $T$ around $P$ which do not fix $Q$ and $x_1,x_3,\ldots$ act as rotations around $Q$ which do not fix $P$, because $a$ is the only elements of the normal form which in $G_e$. Let $d_T$ denote the graph distance in $T$. Then $d_T(Q,a(Q))=0$, $d_T(Q,x_na(Q))=2$, $d_T(Q,x_{n-1}x_na(Q))=2$, $d_T(Q,x_{n-2}x_{n-1}x_na(Q))=4$, $d_T(Q,x_{n-3}x_{n-2}x_{n-1}x_n a(Q))=4$ etc. The case $x_n\in T_2$ is similar. Hence $g$ is not in the kernel of $ \varphi$ and so $ \varphi$ is injective.

Now we assume that $X/G$ is a loop. Since $G$ acts without inversion, there is a $G$-invariant orientation of the edges. Let $e=(P,Q)$ be such an oriented edge and let $t$ be any element of $G$ such that $t(P)=Q$. 
Then $G_{t(e)}t=tG_{e}=\{g\in G\mid g(e)=t(e)\}$. The map $\varphi: G_e\to G_{t(e)}$ given by $\varphi (g)=tgt^{-1}$ is a group isomorphism. 

Let $G'$ denote the HNN-extension $\left<G_Q,t\mid tgt^{-1}=\varphi(g), g\in G_e\right>$. By identifying $G_Q$, $G_e$ and $t$ in $G'$ and $G$ we obtain a unique homomorphism $\psi:G'\to G$. This follows from the fact that any of the relations in $G'$ also hold in $G$. This homomorphism yields an action of $G'$ on $T$.
The action is transitive, because $G_Q\cup \{t\}$ generates $G$. It follows that $\psi$ is surjective.

To show that $\psi$ is injective we choose a $g$ in $G'\setminus \{1\}$. If $g$ is in $G_e$ then $g$ is not in the kernel of $\psi$, because $\psi$ is the identity on $G_e$. If $g$ is $G'\setminus G_e$, then $g$ has a normal form $(x_0,t^{\varepsilon_0},x_1,t^{\varepsilon_1},\ldots ,x_n,t^{\varepsilon_n},h)$ whose length is at least 2, $h$ is in $G_Q$, $x_i$ in $G_e\cup G_{t(e)}$  and $\varepsilon_i\in\{-1,1\}$. Given a vertex $v$ of $T$, then the action of $x_i$ or $h$ does not change the distance to $Q$, because $G_e\cup G_{t(e)}\subset G_Q$. The left multiplication with elements $t^{\varepsilon_n}$ will increase the distance to $Q$ except for the situation where we have subsequence of the form $t^{-1},1,t$ or $t,1,t^{-1}$ in the normal form, which is not possible. This shows that the action of $g$ does not fix $Q$ and so $g$ is not in the kernel of $\psi$. Hence $\psi$ is injective.
\end{proof}

\begin{coro}\label{coro:serre}
A group which acts transitively on an infinite tree splits over the pointwise stabilizer of an edge.
\end{coro}

The barycentric subdivision $T'$ of a graph $T$ is obtained by replacing each edge by a path of length two. In other words, we add an additional vertex on each edge.

\begin{proof}[Proof of Corollary~\ref{coro:serre}]
If there is no edge inversion then $T/G$ is a loop and $G$ splits as HNN-extension according to Theorem~\ref{theo:serre}. If the is an edge inversion then the action of $G$ on the barycentric subdivision $T'$ of $T$ has no edge inversion, the quotient $T/G$ is a segment and $G$ splits as free product with amalgamation according to Theorem~\ref{theo:serre}.
\end{proof}

\section{Ends of groups and Stallings' Theorem}\label{sect:stallings}

Let a group $G$ be generated by $S$. Then $X=\Cay(G,S)$ is defined by $VX=G$ and vertices (group elements) $x,y$ are adjacent if $x^{-1}y\in S$. Equivalently, if there is an $s\in S$ such that $xs=y$. 
A group acts on its Cayley graphs freely by left multiplication. That is, if a group element fixes some vertex then it is the neutral element and hence it fixes all vertices. This implies that the stabilizers of finite sets of vertices are finite subgroups. Right multiplication is only an action a Cayley graph if the group is Abelian.

Given two finite generating sets $S$ and $S'$, the Cayley graphs $\Cay(G,S)$ and $\Cay(G,S')$ are quasi-isometric. That is, if $d_X,d_Y$ denote the graph metrics of $X=\Cay(G,S)$ and $Y=\Cay(G,S')$, respectively, then the identity map $\alpha:G\to G$ is a quasi-isometry between the metric spaces $(G,d_X)$ and $(G,d_Y)$. Formally, there is an integer $a$ such that
\[d_X(x,y)/a-a\le d_Y(\alpha(x),\alpha(y))=d_Y(x,y)\le a d_X(x,y)+a,\]
for all $x,y\in G$. The number of ends in locally finite graphs is a quasi-isometry invariant and hence it does not depend on the finite  generating set. We can speak of the \emph{number of ends of a group}. Finite groups have no ends, infinite finitely generated groups have either one, two or infinitely many ends.

The following criterion is purely algebraic and does not use graphs: An infinite finitely generated group $G$ has more than one end if and only if there is a subset $C$ of $G$ such that $C$ and $G\setminus C$ are infinite and $Cg\setminus C$ is finite, for all $g\in G$ (equivalently, for all $g$ in some set of generators).

\begin{theo}[Stalling's Theorem]
A finitely generated group has more than one end if and only if it splits over some finite subgroup.
\end{theo}

\begin{proof}
Consider a finitely generated group $G$ which splits over some finite subgroup $A$. In the case of an amalgamated product $G=H*_{A}J$ let $S$ be a finite generating set which is contained in $H\cup J$. In the case of an HNN-extension $G=H*^{A}$, let $S$ be a finite generating set of $H$ together with $t$. If $G=H*_{A}J$ then every path from the set $C$ of vertices whose normal presentation starts with an $x_0\in T_H\setminus \{1\}$ is separated by the finite set $A$ from any vertex whose normal presentation starts with an element $x_0\in T_J\setminus \{1\}$, which is the set $G\setminus (C\cup A)$. Using the algebraic definition above, $C$ and $G\setminus C$ are infinite and $Cs\setminus C$ is finite, for all $s\in S$. Hence $\Cay(G,S)$ has more than one end. 

If $G=H*^A$ then every path from the set of vertices $C$ whose normal presentation starts with $t$ is separated from $G\setminus (C\cup A)$ by $A$. Again $\Cay(G,S)$ has more than one end. Hence an infinite finitely generated group which splits over some finite subgroup has more than one end.

Let $X=\Cay(G,S)$ be the Cayley graph with respect to some finite generating set $S$ of a group with more than one end. By Theorem~\ref{theo:optimally} the set $\cc$ of optimally nested cuts is nested and $G$-invariant, because it is invariant under any automorphism. We could as well choose any set $\co$ of optimally nested cuts and set $\cc=\{g(C),g(C^\mathrm{c})\mid C\in\co, g\in G\}$. The transitive action of $G$ on $X$ by left multiplication induces a transient action on $T(\cc)$. By Corollary~\ref{coro:serre}, $G$ splits over a stabilizer of an edge of $T$. Stabilizers of edges in $T$ are stabilizers of edge boundaries $\delta C$ in $X$, for $C\in\cc$. These stabilizers are finite, because the action on $X$ is free.
\end{proof}

\section{Comparison to other papers concerning structure trees}\label{sect:comparison}

In previous papers (for instance \cite{Dunwoody1979,Kroen2001,Moeller1992ends1,Moeller1992ends2,Thomassen1993}), vertices of the structure tree where not defined as inseparable blocks but as equivalence classes of cuts. Cuts $C,D\in\cc$ are called \emph{equivalent} if either (a) $C=D$ or if (b) $C^\mathrm{c}\subset D$ and $C^\mathrm{c}\subset E\subset D$ implies $C^\mathrm{c}=E$ or $E= D$. To prove transitivity of this relation is a bit technical.

It is common to consider so-called structure maps $ \varphi:VX\to VT$ and, for locally finite graphs, $\Phi:\Omega X\to VT\cup \Omega T$. A vertex $x\in VX$ is mapped to the vertex (equivalence class) $v\in VT$ if $x$ is contained in all cuts of this equivalence class. Here it may happen that $ \varphi^{-1}(v)=\emptyset$. Hence there may be vertices of the tree which do not correspond to sets of vertices in the graph. Blocks which correspond to vertices $v$ of the tree where called \emph{regions of $v$} in \cite{Kroen2001}.

\end{document}